\documentclass[12pt]{article}
\usepackage{a4}
\begin{document}
\newtheorem{theorem}{Theorem}
\newcommand{\rk}[1]{\ppar{\bf #1}\quad} 
\newcommand{\ppar}{\par\vskip 8pt plus4pt minus4pt} 
\newcommand{\sh}[1]{\ppar{\bf #1}\par\medskip\nobreak}
\newcommand{\np}{\newpage}            
\newcommand{\nl}{\hfil\break}         
\newcommand{\cl}{\centerline}         
\newcommand{\fnote}[1]{\footnote{\small #1}}


\def\tr{\hbox{tr\,}}
\def\ds{\displaystyle}

\cl{\Large The Burau matrix and Fiedler's invariant for a closed
braid}                

\vglue 0.15truein

\cl{\sc H.R.Morton \fnote{This work was carried out during the low-dimensional
topology program at MSRI, Berkeley in 1996-97. I am grateful to MSRI for their
support.}}

\vglue 0.1truein

{\small\sl 
\cl{Department of Mathematical Sciences,}           
\cl{University of Liverpool,}
\cl{Liverpool L69 3BX, England.}
}

\vglue 0.1truein

\sh{Abstract} 

It is shown how Fiedler's `small state-sum' invariant for a braid $\beta$ can
be calculated from the 2-variable Alexander polynomial of the link which
consists of the closed braid  $\hat{\beta}$  together with the braid axis
$A$.

\rk{AMS Classification numbers}Primary:\quad 57M25. 


\rk{Keywords:} small state-sum, Fiedler's invariant, closed braid, Burau matrix,
2-variable Alexander polynomial.

\bigskip
Original version April 1997, revised November 1997.

\np

\section{Introduction} 

In a recent paper \cite{Fiedler} Fiedler introduced a
simple invariant for a knot $K$ in a line bundle over a surface $F$  by means
of a `small' state-sum, which keeps a count of features of the links resulting
from smoothing each crossing of the projection of $K$ on $F$. The invariant
takes values in a quotient of the integer group ring of $H_1(F)$. Fiedler gives
a number of applications of his general construction. In particular, where $K$
is a closed braid, and can thus be regarded as a knot in a solid torus $V$, his
method gives an invariant of a braid $\beta\in B_n$ in ${\bf Z}[H_1(V)]={\bf
Z}[x^{\pm 1}]$ modulo the relation $x^n=1$. This invariant depends only on the
closure of the braid in $V$ and hence gives an invariant of $\beta$ up to
conjugacy in $B_n$. Its behaviour under Birman and Menasco's
exchange moves has been used to help in detecting when two braids may be
related by such a move.

The purpose of this paper is to show how Fiedler's invariant for a closed braid
$\hat{\beta}$ can be found in terms of the Burau representation of $\beta$, and
hence from the
$2$-variable  Alexander polynomial of the link
$\hat{\beta}\cup A$ consisting of the closed braid $\hat{\beta}$ and its axis
$A$. Its construction here from the Alexander polynomial can be compared with
methods which yield Vassiliev invariants of degree 1 in other contexts, and
suggests possible interpretations of Fiedler's invariants as Vassiliev
invariants of degree 1 in the line bundle.

Having seen how the special case of Fiedler's invariant is related to an
Alexander polynomial I finish the paper with a suggestion of extracting 
similar invariants from the 2-variable Alexander polynomial of a more
general 2-component link. These might be regarded as degree 1 Vassiliev
invariants of one component of the link  when considered as a knot in the
complement of the other component. It would be interesting to know if there was
any similar state sum interpretation of these invariants in the more general
setting.

\section{Burau matrices}

 I make use of the fact that the $2$-variable Alexander
polynomial $\Delta_{\hat{\beta}\cup A}(t,x)$ of a closed braid and its axis can
be calculated as the characteristic polynomial, $\det(I-x\overline{B}(t))$, of
the reduced $(n-1)\times(n-1)$ Burau matrix $\overline{B}(t)$ of the braid
$\beta$, \cite{Morton}.  Since the full $n\times n$ Burau matrix $B(t)$ is
conjugate to $\pmatrix{\overline{B}(t)&\bf{v}\cr\bf{0}&1\cr}$ we can write
$$(1-x)\Delta_{\hat{\beta}\cup A}(t,x)=\det(I-x{B}(t)).$$

Put $t=e^h$ in $\det(I-xB(t))=1+b_1(t)x+\cdots+b_n(t)x^n$, and expand this as a
power series in $h$ to give
$$\det(I-xB(e^h))=\sum_{i=0}^\infty a_i(x)h^i,$$
 where each coefficient $a_i(x)$ is a polynomial in $x$ of degree at most $n$. 

When we set $h=0$, and thus $t=1$, we must get $\Delta_A(x)\times (1-x^n)$ by
the Torres-Fox formula, since the two components $A$ and $\hat{\beta}$ of the
link have linking number $n$. Hence $a_0(x)=1-x^n$. Setting $x=0$ shows also
that $a_1(x)=f_1x+f_2x^2+\cdots+f_nx^n$ for some integers $f_1,\ldots,f_n$. We
know that the determinant of the Burau matrix is $(-t)^{w(\beta)}$, where
$w(\beta)$ is the writhe of the braid, and so $b_n(t)=(-1)^n(-t)^{w(\beta)}$.
Now $w(\beta)=n-1 \bmod 2$ since $\beta$ closes to a single component. Hence
$b_n(e^h)=-1-w(\beta)h+O(h^2)$, giving  $f_n=-w(\beta)$. We shall relate the
remaining coefficients $f_1,\ldots,f_{n-1}$ directly to Fiedler's invariant.

\section{Fiedler's braid invariant.}

 Fiedler's invariant $F_{\beta}$ for an $n$-braid $\beta$ which closes to a
single curve is a symmetric Laurent polynomial, which is even or odd depending
on the parity of $n$. Suppose that the braid \[
\beta=\prod_{r=1}^k\sigma_{i_r}^{\varepsilon_r} \] has been given in terms of
the Artin generators $\sigma_i$, where $\varepsilon_r=\pm1$. Suppose that the
product reads from top to bottom in the braid and the strings are oriented
downwards. For the $r$th crossing define a positive integer $m(r)$ by smoothing
the crossing and following the `ascending string' at the smoothed crossing
around the closed braid until it closes again after $m(r)$ turns around the 
axis. Here the ascending string means the string which starts from the end of
the overcrossing, and is thus string $i_r$ for a positive crossing and string
$i_r +1$ for a negative crossing.  Fiedler's polynomial $F_\beta(X)$ is defined
as a   sum over the $k$ crossings of $\beta$ by  $$F_\beta(X)=\sum_{r=1}^k
\varepsilon_rX^{2m(r)-n}.$$ For a given $m$ we can then write the coefficient
of $X^{2m-n}$ as  $\ds\sum_{m(r)=m} \varepsilon_r$.

\begin{theorem} Let the $n$-string braid $\beta$ have Burau matrix $B(t)$, and
write $\det(I-xB(e^h))= a_0(x)+a_1(x)h+O(h^2)$. Fiedler's polynomial for
$\beta$ satisfies
$$F_\beta(x^{1/2})=(f_1x+\cdots+f_{n-1}x^{n-1})x^{-({n\over 2})},$$ where
$a_1(x)= f_1x+\cdots+f_{n-1}x^{n-1}+f_nx^n$.
\end{theorem}

\noindent{\sl Proof:} Use the classical trace formula for the characteristic
polynomial of a matrix $B$. Suppose that $B$ has eigenvalues
$\lambda_1,\ldots,\lambda_n$. Then $B^m$ has eigenvalues
$\lambda_1^m,\ldots,\lambda_n^m$ and $\det(I-xB)=\prod_{i=1}^n(1-x\lambda_i)$.
Hence \[
\begin{array}{ll}
\ln(\det(I-xB))&\ds=\sum_{i=1}^n
\ln(I-x\lambda_i)=-\sum_{m=1}^\infty\sum_{i=1}^n{1\over m}x^m\lambda_i^m \\ &\ds
=-\sum_{m=1}^\infty {x^m\over m}\tr (B^m),\\
\end{array}\]
 as power series in $x$.

Now expand $\ln(a_0(x)+a_1(x)h+\cdots)$ as a power series in $h$, only  as far
as the term in $h$. We get \[
\begin{array}{l}
\begin{array}{ll}
\ln(a_0(x)+a_1(x)h+\cdots)&=\ds \ln a_0(x)+\ln(1+{a_1(x)\over a_0(x)}h+O(h^2))\\
&
\ds=\ln a_0(x)+{a_1(x)\over a_0(x)}h+O(h^2)\\
\end{array}\\
 =-x^n-x^{2n}/2-\cdots 
+h(f_1x+f_2x^2+\cdots+f_nx^n)(1+x^n+x^{2n}+\cdots)+O(h^2).\\
\end{array}\]

The trace formula above applied to $B(e^h)$ shows at once that
$\tr((B(e^h))^m)=-mf_mh+O(h^2)$ for $1\le m<n$.

 The proof will be completed by relating the term in $h$  in the trace of this
matrix to the appropriate coefficient of Fiedler's polynomial. It is thus 
enough to show that $\ds\tr((B(e^h))^m)=-m(\sum_{m(r)=m}\varepsilon_r)h+O(h^2)$
for $1\le m<n$.

The Burau representation $\rho:B_n\to GL(n,{\bf Z}[t^{\pm 1}])$ is the group
homomorphism defined on the elementary braid $\sigma_i$ by
\[
\rho(\sigma_i)=B_i=\pmatrix{\begin{array}{c| c| c}I_{i-1}&0&0\\ \hline
 0&\begin{array}{cc}1-t& t\cr 1&0\cr\end{array}&0\cr\hline
0&0&I_{n-i-1}\cr\end{array}}.
\] The Burau matrix for the given braid $\beta$ is then
\[ B(t)=\rho(\beta)=\prod_{r=1}^k B_{i_r}^{\varepsilon_r}.
\]

Now \[\begin{array}{ll}
B_i(e^h)&=\pmatrix{\begin{array}{c|c|c}I_{i-1}&0&0\cr\hline
 0&\begin{array}{cc}0& 1\cr 1&0\cr\end{array}&0\cr\hline
0&0&I_{n-i-1}\cr\end{array}} +h\pmatrix{\begin{array}{c|c|c}
 0_{i-1}&0&0\cr\hline
 0&\begin{array}{rr}-1&1\cr0&0\cr\end{array}&0\cr\hline 0&0&0_{n-i-1}\cr
\end{array} }+O(h^2)\\ &=T_i+hP_i^++O(h^2), \hbox{ say}.\\
\end{array}
\] We can similarly write $B_i^{-1}=T_i+hP_i^- +O(h^2)$ where 
\[ P^-_i=\pmatrix{\begin{array}{c|c|c}
 0_{i-1}&0&0\cr\hline
 0&\begin{array}{rr}0&0\cr-1& 1\cr\end{array}&0\cr\hline 0&0&0_{n-i-1}\cr
\end{array} }.
\] Then \[( B(e^h))^m=(\prod_{r=1}^k(T_{i_r}+hP_{i_r}^\pm))^m+O(h^2).
\]

We can write a matrix of the form $M=\prod_{r=1}^l (C_r+hD_r)$ as \[
M=C_1C_2\ldots C_l+h(D_1C_2\ldots C_l+C_1D_2C_3\ldots C_l+\cdots+C_1C_2\ldots
C_{l-1}D_l)+O(h^2),\] and then \[ \tr M=\tr( C_1C_2\ldots
C_l)+h(\tr(D_1C_2\ldots C_l)+\tr(C_1D_2C_3\ldots C_l)+\cdots)+O(h^2).
\] The term in $h$ can be rewritten as
\[\tr(C_2\ldots C_lD_1)+\tr(C_3\ldots C_lC_1D_2)+\cdots+\tr(C_1C_2\ldots
C_{l-1}D_l)\] by cycling the matrices so that the $r$th product  ends with the
matrix $D_r$.

 Apply this to find the term in $h$ in $\tr((B(e^h))^m$ as the sum of $mk$
terms, each of which is the trace of the product of $mk$ matrices of the form
$T_{i_{r+1}}\ldots T_{i_{r-1}}P^\pm_{i_r}$ with sign $\pm$ according to the
sign of $\varepsilon_r$. For each of the $k$ crossings  of the original braid
the matrix $T_{i_{r+1}}\ldots T_{i_{r-1}}P^\pm_{i_r}$ occurs $m$ times in the
sum. Thus
 \[  f_m=-\sum_{r=1}^k\tr(T_{i_{r+1}}\ldots T_{i_{r-1}}P^\pm_{i_r}).
\]

The proof of theorem 1 will be completed by showing that \[
tr(T_{i_{r+1}}\ldots T_{i_{r-1}}P^\pm_{i_r})=\cases{-\varepsilon_r& if
$m(r)=m$\cr 0& otherwise.}\]

The matrix $T_i$ is the permutation matrix for the transposition $(i\;i+1)$.
Hence a product of these matrices is also a permutation matrix, $T$ say, whose
permutation is the product $\pi$ of the corresponding transpositions. Then the
entries in $T$ satisfy 
\[ T_{ij}=\cases{1& if $i=\pi(j)$,\cr 0& otherwise.}\]

The matrix $T=T_{i_{r+1}}\ldots T_{i_{r-1}}$ above is thus a permutation matrix
with permutation $\pi_r^{(m)}$, say. Notice that the permutation corresponding
to the product $TT_{i_r}$ is conjugate to the permutation of the braid
$\beta^m$. Under the assumption that $\beta$ closes to a single curve this will
be the $m$th power of an $n$-cycle, and will hence not fix any number when
$1\le m<n$.  Hence  $\pi_r^{(m)}$ can not carry $i_r$ to $i_r+1$ or vice versa,
in this range.

If the $r$th crossing is smoothed and the strings $i_r$ and $i_r+1$ are
followed upwards around the braid $m$ times, with $m<n$, they will not pass
through the smoothed crossing. They then become the strings $\pi_r^{(m)}(i_r)$
and $\pi_r^{(m)}(i_r+1)$ respectively when they return to the level of the
bottom of the $r$th crossing.  Now when $\varepsilon_r=+1$ the ascending string
at the $r$th crossing, which is string $i_r$,  returns to position $i_r$ after
the permutation $\pi_r^{(m)} $
 if and only if $m=m(r)$. Similarly when $\varepsilon_r=-1$ the ascending
string, in this case string $i_r+1$ returns to position $i_r+1$ exactly when
$m=m(r)$. 

 The matrices $P_{i_r}^\pm$ have only two non-zero entries. Suppose first that
$\varepsilon_r=+1$. Then $\tr(TP_{i_r}^+)$ is the sum of two terms. The
off-diagonal entry gives a contribution only if the permutation matrix $T$ maps
it onto the diagonal. This requires $\pi_r^{(m)}(i_r)=i_r+1$, which was
excluded above. The diagonal entry contributes $-1$ if and only 
$\pi_r^{(m)}(i_r)=i_r$, which is the condition that $m=m(r)$. Thus when
$\varepsilon_r=+1$ we get a contribution of $-\varepsilon_r$ to  the trace if
and only if $m=m(r)$, and zero otherwise. 

A similar argument holds when $\varepsilon_r=-1$. Again the off-diagonal entry
does not contribute to the trace, while the diagonal entry contributes $+1$  if
and only if $\pi_r^{(m)}(i_r+1)=i_r+1$. This corresponds once more to the
condition that $m=m(r)$, and so in each case we have a contribution of
$-\varepsilon_r$ if and only if $m=m(r)$. The total coefficient of $h$ in
$\tr((B(e^h)^m)$ is then $\ds -m\sum_{m=m(r)}\varepsilon_r$, showing that $\ds
f_m=\sum_{m=m(r)}\varepsilon_r$ as claimed. This completes the proof of theorem
1.

\section{Determination from an Alexander polynomial.}

 If we are given the
Alexander polynomial of the closed braid $\hat{\beta}$ and its axis $A$ as a
2-variable polynomial we can recover Fiedler's invariant for the braid. First
multiply by $1-x$, where $x$ is the variable for the axis.  This gives the
characteristic polynomial of the Burau matrix for $\beta$, up to multiplication
by a power of $ x$ and a power of $t$, and a sign. Put $t=e^h$ and expand as a
power series in $h$ with coefficients depending on $x$. Then multiply by a
power of $x$ and a sign to make the constant term $1-x^n$. The result will be
the characteristic polynomial used above, up to a power of $t=e^h$. Extract the
coefficient $f_0+f_1x+\cdots+f_{n-1}x^{n-1}+f_nx^n$ of $h$. This will contain
the Fiedler polynomial as before in the terms $f_1x+\cdots+f_{n-1}x^{n-1}$,
while the remaining terms will come from a factor of $t^{f_0}$ and will satisfy
$f_0+f_n=-w(\beta)$.

A similar interpretation looks plausible for  the coefficients of the linear
terms in $h_1,\ldots,h_k$ when the Alexander polynomial of a closed braid with
$k$ components and its axis is expanded in terms of the meridian generator $x$
for the axis and meridians $t_i=e^{h_i}$ for the components. This polynomial
can again be written in terms of the characteristic polynomial of a suitable
`coloured' Burau matrix. The eventual coefficient of $h_i$ should then have
contributions from the overcrossings of the corresponding component of the
closed braid, as in the Fiedler polynomial above.

As a possible extension to the case of a general link $L$ with two components
$X$ and $T$ say, we might put $t=e^h$ in the Alexander polynomial $\Delta_{X\cup
T}(x,t)$ of $L$ and consider only the terms $a_0(x)+a_1(x)\,h$ up to degree 1
in $h$. The polynomial $a_0(x)$ is $\Delta_X(x)(1-x^n)/(1-x)$, where $n$ is the
linking number of $X$ and $T$, and $\Delta_X(x)$ is the Alexander polynomial of
$X$.  Now consider $a_1(x)$ as a polynomial modulo the ideal generated by
$a_0(x)$. This is an invariant of $L$ as it is unaffected by any ambiguity of
powers of $x$ and $t$ in the Alexander polynomial. This seems to me to be the
nearest analogue to Fiedler's invariant for the link component $T$ with
meridian $t$ when regarded as a knot in the complement of $X$; in the case of a
closed braid we take $X$ as the braid axis and $T$ as the closed braid. It
looks likely to be a Vassiliev invariant of type 1 for knots in the complement
of $X$. There is not, however, any obvious candidate for a state-sum
construction of this invariant along Fiedler's lines when the component $X$ is
knotted.

\end{document}